\documentclass[12pt,twoside]{article}

\makeatother

\normalsize
\usepackage{cite}
\usepackage{amsfonts}
\usepackage[mathscr]{euscript}

\newcommand{\lbl}[1]{\label{#1}}

\newtheorem{re}{Remark}[section]

\newtheorem{theo}{Theorem}[section]
\newtheorem{lem}{Lemma}[section]
\newtheorem{col}{Corollary}[section]

\newcommand{\be}{\begin{equation}}
\newcommand{\ee}{\end{equation}}
\newcommand\bes{\begin{eqnarray}} \newcommand\ees{\end{eqnarray}}
\newcommand{\bess}{\begin{eqnarray*}}
\newcommand{\eess}{\end{eqnarray*}}
\newcommand{\D}{\displaystyle}

\def\theequation{\arabic{section}.\arabic{equation}}

\begin{document}
\setlength{\baselineskip}{16pt} \pagestyle{myheadings}

\begin{center}{\Large\bf  Regularity  criteria for the 3D  MHD equations\\[2mm] in term of velocity}\footnote{Supported by the NNSF of China (nos.\,11271379, 11271381) and the National Basic Research Program of China (973 Program) (Grant No.\,2010CB808002).}\\[2mm]

 \large  Qunyi Bie$^{1,2}$, ~~Qiru Wang$^{1,}$\footnote{Corresponding author. E-mail address: mcswqr@mail.sysu.edu.cn, tel.: +86-20-84037100, fax: +86-20-84037978.},~~Zhengan Yao$^{1}$\\[1mm]

{\small 1. School of Mathematics and Computational Science,  Sun Yat-Sen University,\\
Guangzhou 510275, PR China\\2. College of Science, China Three Gorges University, Yichang 443002, PR China
}\\[1mm]
 \end{center}

 \begin{quote}
 \noindent {\bf Abstract.} {\small In  this  paper we consider three-dimensional incompressible
magnetohydrodynamics equations. By using interpolation inequalities in anisotropic Lebesgue space, we provide regularity criteria involving  the velocity  or alternatively involving the fractional derivative of velocity in one direction, which generalize  some known results.}\\[2mm]
 \noindent {\bf Key wor{\rm d}s:} 3D Magnetohydrodynamics equations; regularity criteria; anisotropic Lebesgue space\\[2mm]

 \end{quote}
 \setcounter{equation}{0}
 \def\theequation{\arabic{section}.\arabic{equation}}
 \vskip 4pt
 \section{Introduction and main results} \setcounter{equation}{0}

In this paper, we consider the three dimensional  magnetohydrodynamics (MHD) equations:
 \bes
 \left\{\begin{array}{ll}\displaystyle \partial_tu+(u\cdot\nabla) u
-\nu\Delta u+\nabla p=(b\cdot\nabla)b, &\ \ \ {\rm
in}\ \ \ \mathbb{R}^3\times(0,\infty),\\
 \medskip \displaystyle \partial_t b +(u\cdot\nabla)b-\eta\Delta b =(b\cdot\nabla) u, &\ \ \ {\rm
in}\ \ \mathbb{R}^3\times(0,\infty),\\
 {\rm div} u=0,\ \ {\rm div} b=0, &\ \ \ {\rm in}\ \
 \mathbb{R}^3\times(0,\infty),\\
u(x,0)=u_0(x),~b(x,0)=b_0(x),&\ \ \ {\rm in}\ \ \mathbb{R}^3,
 \end{array}
 \right.
 \label{1.1}
 \ees
where $u=(u_1, u_2, u_3), b=(b_1, b_2, b_3)$  and $p$ are  velocity field,
the magnetic field and  the scalar  pressure, respectively.
The parameters $\nu>0, \eta>0$ denote the kinematic viscosity and the magnetic diffusivity, respectively. The incompressible MHD equations describe the motion of an electrical conducting fluid in the presence of a magnetic field and have an important meaning in physics and other applied
areas such as geophysics, astrophysics, and engineering problems (see \cite{L1}).

The local well-posedness of the Cauchy problem (\ref{1.1}) in the usual Sobolev spaces
$H^s(\mathbb{R}^3)$ was established in \cite{ST} for any given initial data $(u_0, b_0)\in H^s(\mathbb{R}^3), s\geq 3$. However, global regularity is still an open problem. There were numerous important progresses on its fundamental issue of the regularity
for the weak solution to (\ref{1.1}) (see\cite{Z1,CW, JL, CMZ, HX1, HX2, G, HW1, HW2, Z6, W1, W2, W3}). We would like to mention that  He and Xin  \cite{HX1, HX2} realized that the velocity fields played a dominate role in the regularity of the solutions to 3D incompressible MHD equations,  and  established the global regularity of the strong solution involving only the velocity field for the first time.

Very recently,  Jia and Zhou \cite{JZ} established the global regularity criterion for (\ref{1.1}) as follows:
 \bes
u_3, b_i\in L^\beta(0, T; L^\alpha(\mathbb{R}^3))~(i=1,2,3),
~~~{\rm with~~} \frac{2}{\beta}+\frac{3}{\alpha}\leq \frac{3}{4}+\frac{1}{2\alpha},~~
\frac{10}{3}<\alpha\leq \infty.
 \lbl{1.10}
 \ees
And in \cite{LD}, Lin and Du showed that the weak solution remains smooth on $(0, T]\times\mathbb{R}^3$ if the derivatives of the velocity in one direction satisfy the condition
 \bes
\int_0^T{\Big\|}\frac{\partial u(s)}{\partial x_i}{\Big\|}_{L^\alpha}^\beta {\rm d}s<\infty, ~~~{\rm with~~}\alpha>2,~~ \frac{2}{\beta}+\frac{3}{\alpha}\leq\frac{3}{4}+\frac{3}{2\alpha},
 \lbl{1.11}
 \ees
for some $i=1,2,3$ and $T>0$.

 The purpose of this paper is to generalize the recent results in \cite{JZ,LD} and establish two sufficient conditions for the global regularity of strong solutions to system (\ref{1.1}). Our idea and  proof framework  of main results come from \cite{CT1, CT2, Z}, in which the authors established some important regularity criteria for the 3D Navier-Stokes equations. Before stating our main results, we introduce some basic functional spaces and the definitions of the weak and strong solutions.

We denote by $L^q$ and $H^m$ the usual $L^q$-Lebesgue and Sobolev spaces, respectively. Set
$$\mathcal{V}=\{\phi: {\rm ~the~ 3D ~ vector~ valued} ~C_0^\infty~{\rm functions~and~}\nabla\cdot \phi=0 \},$$
which will compose the space of test functions. Let $H$ and $V$ be the closure spaces of $\mathcal{V}$ in $L^2$ under $L^2$-topology, and in $H^1$ under  $H^1$-topology, respectively.
A pair $(u,b)$ of measurable functions is called {\it a weak solution} to {\rm(\ref{1.1})}
with $(u_0, b_0)\in H$, provided that $(u, b)$ satisfies the following three conditions:

{\rm (1)} $(u,b)\in C_w([0,T], H)\cap L^2(0, T; V)$, ~{\rm and~}$(\partial_t u, \partial_t b)\in L^1(0, T;  V^\prime)$,~where $V^\prime$ is the dual space of $V$;

{\rm (2)} the weak formulations of the  MHD equations:
 \bess
 \begin{array}{ll}
&\hspace{-0.5cm}\displaystyle\int_{\mathbb{R}^3}u(x,t)\varphi(x,t){\rm d}x-\int_{\mathbb{R}^3}u(x,t_0)\varphi(x,t_0){\rm d}x\\[3mm]
&=\displaystyle\int_{t_0}^t\int_{\mathbb{R}^3}u(x,s)(\varphi_t(x,s)+\nu\Delta \varphi (x,s)){\rm d}x{\rm d}s \\[3mm]
&\ \ \  +\displaystyle\int_{t_0}^t\int_{\mathbb{R}^3}(u(x,s)\cdot\nabla)\varphi(x,s)\cdot u(x,s){\rm d}x{\rm d}s\\[3mm]
&\ \ \  -\displaystyle\int_{t_0}^t\int_{\mathbb{R}^3}(b(x,s)\cdot\nabla)\varphi(x,s)\cdot b(x,s){\rm d}x{\rm d}s,
  \end{array}
 \eess

 \bess
 \begin{array}{ll}
&\hspace{-0.5cm}\displaystyle\int_{\mathbb{R}^3}b(x,t)\varphi(x,t){\rm d}x-\int_{\mathbb{R}^3}b(x,t_0)\varphi(x,t_0){\rm d}x\\[3mm]
&=\displaystyle\int_{t_0}^t\int_{\mathbb{R}^3}b(x,s)(\varphi_t(x,s)+\eta\Delta \varphi (x,s)){\rm d}x{\rm d}s \\[3mm]
&\ \ \ + \displaystyle\int_{t_0}^t\int_{\mathbb{R}^3}(u(x,s)\cdot\nabla)\varphi(x,s)\cdot b(x,s){\rm d}x{\rm d}s\\[3mm]
&\ \ \ -\displaystyle\int_{t_0}^t\int_{\mathbb{R}^3}(b(x,s)\cdot\nabla)\varphi(x,s)\cdot u(x,s){\rm d}x{\rm d}s
 \end{array}
 \eess
for almost every $t,t_0\in [0, T]$ and for every test function $\varphi \in C^\infty([0, T], \mathcal{V})$;

{\rm (3)} the energy inequality:
 \bes
 \begin{array}[b]{ll}
&\hspace{-0.5cm}\displaystyle\|u(t)\|_{L^2}^2+\|b(t)\|_{L^2}^2+\nu\int_{t_0}^t\|\nabla u(\tau)\|_{L^2}^2{\rm d}\tau
+\eta\int_{t_0}^t\|\nabla b(\tau)\|_{L^2}^2{\rm d}\tau\\[3mm]
&\leq \|u(t_0)\|_{L^2}^2+\|b(t_0)\|_{L^2}^2
 \lbl{1.2}
 \end{array}
 \ees
for every $t$ and almost every $t_0$.

In addition, if $(u_0, b_0)\in V$, a weak solution is said to be {\it a strong solution} to {\rm (\ref{1.1})},
provided
$$u\in C([0, T], V)\cap L^2(0, T; H^2), ~~{\rm and}~~~(\partial_t u, ~\partial_t b)\in L^2(0, T; H).$$

Throughout the paper, we denote $\nabla _h=(\partial_1, \partial_2)$
and $(i, j, k)$ belongs to permutation group $S_3$:=span\{1, 2, 3\}.
We also agree that $\Lambda_i=\sqrt{-\partial_i^2}$, ~$\Lambda=\sqrt{-\Delta}$~~and
$${\Big\|}\|f\|_{L_i^p}{\Big\|}_{L_{j,k}^q}
:={\Big(}\int_{\mathbb{R}^2}{\Big(}\int_{\mathbb{R}}|f(x)|^p{\rm d}x_i{\Big)}^{\frac{q}{p}}{\rm d}x_j{\rm d}x_k{\Big)}^{\frac{1}{q}}.$$
For simplicity, we assume $\nu=\eta=1$ throughout this paper. Denote
 \bes
 \left.
 \begin{array}{ll}
 \mathcal{G}_1={\Big\{}(p,\alpha)\in (2,\infty]^2; ~~\displaystyle\frac{1}{p}+\frac{2}{\alpha}<1 ~~{\rm and}~~ \frac{\alpha}{p(\alpha-2)}<\frac{3}{4}{\Big\}},\\[3mm]
 \mathcal{G}_2={\Big\{}\alpha\in (2,\infty];
~~\displaystyle\frac{3-\alpha\gamma}{\alpha}<1 ~~{\rm and}~~ \frac{1-\alpha\gamma}{\alpha-2}<\frac{3}{4} {\Big\}}.
 \end{array}
 \right.
 \nonumber
 \ees

The main results in this paper are stated as follows.
 \begin{theo}
 \lbl{th1.1}
Let $(u_0, b_0)\in H^3$ with ${\rm div} u_0={\rm div} b_0=0$, $(u, b)$
be a weak solution of {\rm (\ref{1.1})}. If $u(x, t)$ satisfies the following condition
 \bes
 \int_0^T{\Big\|}\|u(\tau)\|_{L_i^p}{\Big\|}_{L_{j, k}^{\alpha}}^{\beta}{\rm d}\tau<\infty,
~~{\rm with}~~\frac{2}{\beta}
+\frac{2}{\alpha}+\frac{1}{p}\leq \frac{3}{4}+\frac{1}{2\alpha},~~(p, \alpha)\in \mathcal{G}_1,
 \lbl{1.3}
 \ees
 then $(u, b)$ is the strong solution of {\rm (\ref{1.1})} on $[0, T]$.
 \end{theo}
 \begin{re}
When  we fix  $p=\alpha$, compared with  the sufficient condition {\rm (\ref{1.10})} as stated in {\rm \cite{JZ}}, the condition {\rm (\ref{1.3})} involves only the velocity fields.
 \end{re}

As a corollary of Theorem \ref{th1.1}, we can further obtain the following criterion in terms of
the fractional derivative of velocity field in one direction.

 \begin{col}
 \lbl{c1}
Let $(u_0, b_0)\in H^3$ with ${\rm div} u_0={\rm div} b_0=0$, and $(u, b)$
be a weak solution of {\rm (\ref{1.1})}. Suppose that for some $i\in \{1, 2, 3\}$,
 \bes
 \displaystyle\int_0^T\|\Lambda_i^\gamma u(\tau)\|_{L^\alpha(\mathbb{R}^3)}^\beta {\rm d}\tau<\infty,~\left\{\begin{array}{ll}\D \frac{2}{\beta}+\frac{3}{\alpha}\leq \frac{3}{4}
+\frac{3}{2\alpha},~&\alpha\in (2,\infty), \D ~{\rm if}~\gamma\in (\frac{1}{\alpha}, 1];\\[3mm]
 \D\frac{2}{\beta}+\frac{3}{\alpha}\leq \frac{3}{4}+\gamma
+\frac{1}{2\alpha},~&\alpha\in \mathcal{G}_2, \D ~{\rm if}~\gamma\in [0,\frac{1}{\alpha}).
 \end{array}
 \right.
 \lbl{1.4}
 \ees
Then $(u, b)$ is the strong solution of {\rm (\ref{1.1})} on $[0, T]$.
 \end{col}
\begin{re}
If  we take $\gamma=1$, then the sufficient condition {\rm (\ref{1.4})} naturally turn into condition {\rm (\ref{1.11})}. Thus  condition  {\rm (\ref{1.4})} is a generalization of condition {\rm (\ref{1.11})} as stated in {\rm \cite{LD}} in terms of derivatives of the velocity fields in one direction.
 \end{re}

 \section{Proof of main results} \setcounter{equation}{0}
In this section, we give complete proofs of the results described in Section 1. To do this, we adopt the similar proof framework used in \cite{CT1, CT2, Z}. Firstly we recall the following inequalities which may be found in \cite{A, DL, L} (see also \cite{Z}):
 \begin{lem}
 \lbl{lm1.1}
Let $2\leq p\leq 6, ~2\leq q< \infty$ and ~$2<r, s\leq \infty$. There hold that,
 \bes
\|f\|_{L^p(\mathbb{R}^3)}\leq C\|f\|_{L^2(\mathbb{R}^3)}^{\frac{6-p}{2p}}
\|\partial_1f\|_{L^2(\mathbb{R}^3)}^{\frac{p-2}{2p}}
\|\partial_2f\|_{L^2(\mathbb{R}^3)}^{\frac{p-2}{2p}}
\|\partial_3f\|_{L^2(\mathbb{R}^3)}^{\frac{p-2}{2p}},
 \nonumber
 \ees
 \bes
{\Big\|}\|f\|_{L_{i, j}^q}{\Big\|}_{L_{k}^2}\leq C\|f\|_{L^2(\mathbb{R}^3)}^{\frac{2}{q}}
\|\partial_if\|_{L^2(\mathbb{R}^3)}^{\frac{q-2}{2q}}
\|\partial_jf\|_{L^2(\mathbb{R}^3)}^{\frac{q-2}{2q}}
 \nonumber
 \ees
and
 \bes
{\Big\|}\|f\|_{L_i^{\frac{2r}{r-2}}}{\Big\|}_{L_{j, k}^{\frac{2s}{s-2}}}\leq C\|f\|_{L^2(\mathbb{R}^3)}^{1-\frac{1}{r}-\frac{2}{s}}\|\partial_if\|_{L^2(\mathbb{R}^3)}^{\frac{1}{r}}
\|\partial_jf\|_{L^2(\mathbb{R}^3)}^{\frac{1}{s}}\|\partial_kf\|_{L^2(\mathbb{R}^3)}^{\frac{1}{s}}.
 \nonumber
 \ees
 \end{lem}

 \subsection{Proof of Theorem \ref{th1.1}}
Similar to the Navier-Stokes equation, the global existence of weak solutions to the MHD equations can be proved by applying Galerkin's method and compact argument, see \cite{DL}.  It is also well known that there exists a unique strong solution for a short time interval  if $(u_0, b_0)\in V$. Furthermore, this strong solution is the only weak solution with the initial value $(u_0, b_0)$, on  the maximal interval of existence of the strong solution.

Let $(u, b)$ be the strong solution with the initial value $(u_0, b_0)\in V$ such that $(u, b)\in C([0, T^\ast), V)\cap L^2([0, T^\ast), H^2)$, where $[0, T^\ast)$ is the maximal interval of existence of the unique strong solution. When $T^\ast\geq T$, there is nothing to prove, while when $T^\ast < T$, then our strategy is to show that the $H^1$ norm of this strong solution is bounded on $[0, T^\ast)$, provided condition (\ref{1.3}) or (\ref{1.4}) hol{\rm d}s. Thus the interval  $[0, T^\ast)$ can not be a maximal interval of existence, and this lea{\rm d}s to a contradiction.

In the following we assume that $(u, b)$ is the strong solution on its maximal interval of existence $[0, T^\ast)$ satisfying $T^\ast < T$. Recall that the strong solution is indeed the weak solution on  $[0, T^\ast)$. Therefore, by the energy inequality (\ref{1.2}), $(u, b)$ satisfies
 \bes
 \|u(t)\|_{L^2}^2+\|b(t)\|_{L^2}^2+\nu\int_0^t\|\nabla u(\tau)\|_{L^2}^2{\rm d}\tau
+\eta\int_0^t\|\nabla b(\tau)\|_{L^2}^2{\rm d}\tau\leq C,
 \lbl{2.0}
 \ees
where $C=\|u_0\|_{L^2}^2+\|b_0\|_{L^2}^2$.

Without loss of generality, we consider the case
\bes
 \int_0^t{\Big\|}\|u(\tau)\|_{L_3^p}{\Big\|}_{L_{1, 2}^{\alpha}}^{\beta}{\rm d}\tau<\infty,
~~{\rm with}~~\frac{2}{\beta}
+\frac{2}{\alpha}+\frac{1}{p}\leq \frac{3}{4}+\frac{1}{2\alpha},~~(p, \alpha)\in \mathcal{G}_1.
 \lbl{2.30}
 \ees
Next, let us show that the $H^1$ norm of the strong solution $(u, b)$ is bounded on interval $[0, T^\ast)$.
We start with the estimates of $\|\nabla_h u\|_{L^2}$ and $\|\nabla_h b\|_{L^2}$.

{\em Step 1. Estimates for $\|\nabla_h u\|_{L^2}$ and $\|\nabla_h b\|_{L^2}$}.
Taking the inner product of the first equation and second equations in (\ref{1.1})
with $-\Delta_h u$ and $-\Delta_h b$ in $H$, respectively, we have
 \bes
 \begin{array}[b]{ll}
&\hspace{-0.5cm} \D\frac{1}{2}\frac{d}{dt}{\Big(}\|\nabla_h u(t)\|_{L^2}^2
+\|\nabla_h b(t)\|_{L^2}^2{\Big)}+{\Big(}\|\nabla_h \nabla u(t)\|_{L^2}^2+\|\nabla_h\nabla b(t)\|_{L^2}^2{\Big)}\\[3mm]
&=\D\int_{\mathbb{R}^3}(u\cdot\nabla)u\cdot \Delta_h u{\rm d}x-\int_{\mathbb{R}^3}(b\cdot\nabla)b\cdot \Delta_h u{\rm d}x \\[3mm]
&\ \ \ +\D\int_{\mathbb{R}^3}(u\cdot\nabla)b\cdot \Delta_h b{\rm d}x-\int_{\mathbb{R}^3}(b\cdot\nabla)u\cdot \Delta_h b{\rm d}x\\[3mm]
&=:I_1+I_2+I_3+I_4.
 \lbl{2.1}
 \end{array}
 \ees
Attention is now focused on bounding these terms; we start with $I_1$.
By integration by parts and the incompressible conditions, one gets
 \bes
I_1&=&\int_{\mathbb{R}^3}(u\cdot\nabla)u\cdot \Delta_h u{\rm d}x
=\int_{\mathbb{R}^3}\sum_{j=1}^3\sum_{l=1}^2u_j\partial_ju\cdot \partial_l^2 u{\rm d}x\nonumber\\
&=&-\int_{\mathbb{R}^3}\sum_{j=1}^3\sum_{l=1}^2\partial_lu_j\partial_ju\cdot \partial_l u{\rm d}x\nonumber \\
&=&\int_{\mathbb{R}^3}\sum_{j=1}^3\sum_{l=1}^2u_j\partial_l(\partial_ju\cdot \partial_lu){\rm d}x\nonumber\\
&\leq& \int_{\mathbb{R}^3}|u||\nabla u||\nabla_h\nabla u|{\rm d}x,
 \lbl{2.20}
 \ees
where we have used
$$\int_{\mathbb{R}^3}\sum_{j=1}^3\sum_{l=1}^2u_j\partial_j\partial_l u\cdot \partial_l u{\rm d}x=0.$$
Similarly, we can estimate the terms $I_2, I_3$ and $I_4$ as follows.
 \bes
 \begin{array}[b]{ll}
I_2+I_4&=\D -\int_{\mathbb{R}^3}(b\cdot\nabla)b\cdot \Delta_h u{\rm d}x
-\int_{\mathbb{R}^3}(b\cdot\nabla)u\cdot \Delta_h b{\rm d}x\\[3mm]
&=\D-\int_{\mathbb{R}^3}\sum_{j=1}^3\sum_{l=1}^2b_j\partial_jb\cdot \partial_l^2 u{\rm d}x
-\int_{\mathbb{R}^3}\sum_{j=1}^3\sum_{l=1}^2b_j\partial_ju\cdot \partial_l^2 b{\rm d}x\\[3mm]
&=\D\int_{\mathbb{R}^3}\sum_{j=1}^3\sum_{l=1}^2\partial_lb_j\partial_jb\cdot \partial_l u{\rm d}x+\int_{\mathbb{R}^3}\sum_{j=1}^3\sum_{l=1}^2b_j\partial_j\partial_lb\cdot \partial_l u{\rm d}x\\[3mm]
&\ \ \ +\D\int_{\mathbb{R}^3}\sum_{j=1}^3\sum_{l=1}^2\partial_lb_j\partial_ju\cdot \partial_lb{\rm d}x+\int_{\mathbb{R}^3}\sum_{j=1}^3\sum_{l=1}^2b_j\partial_j\partial_lu\cdot \partial_l b{\rm d}x\\[3mm]
&=\D-\int_{\mathbb{R}^3}\sum_{j=1}^3\sum_{l=1}^2u\cdot \partial_l(\partial_lb_j\partial_jb){\rm d}x
-\int_{\mathbb{R}^3}\sum_{j=1}^3\sum_{l=1}^2u\cdot \partial_j(\partial_lb_j\partial_lb){\rm d}x\\[3mm]
&\leq\D \int_{\mathbb{R}^3}|u||\nabla b||\nabla_h\nabla b|{\rm d}x,
 \lbl{2.21}
 \end{array}
 \ees
where we have used
 \bes
&&\hspace{-1cm}\int_{\mathbb{R}^3}\sum_{j=1}^3\sum_{l=1}^2b_j\partial_j\partial_lb\cdot \partial_lu{\rm d}x+\int_{\mathbb{R}^3}\sum_{j=1}^3\sum_{l=1}^2b_j\partial_j\partial_lu\cdot \partial_l b{\rm d}x\nonumber\\[3mm]
&=& \int_{\mathbb{R}^3}\sum_{j=1}^3\sum_{l=1}^2b_j\partial_j(\partial_lb\cdot \partial_l u){\rm d}x=0
 \nonumber
 \ees
due to $\nabla\cdot b=0$. Similar to the estimate for $I_1$, we have
 \bes
I_3&=&\int_{\mathbb{R}^3}(u\cdot\nabla)b\cdot \Delta_h b{\rm d}x
=\int_{\mathbb{R}^3}\sum_{j=1}^3\sum_{l=1}^2u_j\partial_jb\cdot \partial_l^2 b{\rm d}x\nonumber\\[3mm]
&=&-\int_{\mathbb{R}^3}\sum_{j=1}^3\sum_{l=1}^2\partial_lu_j\partial_jb\cdot \partial_l b{\rm d}x\nonumber \\[3mm]
&=&\int_{\mathbb{R}^3}\sum_{j=1}^3\sum_{l=1}^2u_j\partial_l(\partial_jb\cdot \partial_lb){\rm d}x\nonumber\\[3mm]
&\leq& \int_{\mathbb{R}^3}|u||\nabla b||\nabla_h\nabla b|{\rm d}x.
 \lbl{2.22}
 \ees
Combining the above estimates with (\ref{2.1}), it is clear that
 \bes
&&\hspace{-1cm}\frac{1}{2}\frac{d}{dt}{\Big(}\|\nabla_h u(t)\|_{L^2}^2
+\|\nabla_h b(t)\|_{L^2}^2{\Big)}+{\Big(}\|\nabla_h \nabla u(t)\|_{L^2}^2
+\|\nabla_h\nabla b(t)\|_{L^2}^2{\Big)}\nonumber\\[3mm]
&\leq&\int_{\mathbb{R}^3}|u||\nabla u||\nabla_h\nabla u|{\rm d}x
+\int_{\mathbb{R}^3}|u||\nabla b||\nabla_h\nabla b|{\rm d}x\nonumber\\[3mm]
&=:&K_1+K_2.
 \lbl{2.2}
 \ees
Now we estimate $K_1$ and $K_2$. Employing Lemma \ref{lm1.1}, we have that for $(p, \alpha)\in \mathcal{G}_1$,
 \bes
{\Big\|}\|\nabla u\|_{L_3^{\frac{2p}{p-2}}}{\Big\|}_{L_{1, 2}^{\frac{2\alpha}{\alpha-2}}}\leq C\|\nabla u\|_{L^2(\mathbb{R}^3)}^{1-\frac{1}{p}-\frac{2}{\alpha}}\|\partial_1\nabla u\|_{L^2(\mathbb{R}^3)}^{\frac{1}{\alpha}}
\|\partial_2\nabla u\|_{L^2(\mathbb{R}^3)}^{\frac{1}{\alpha}}\|\partial_3\nabla u\|_{L^2(\mathbb{R}^3)}^{\frac{1}{p}}.
 \lbl{2.3}
 \ees
Thus, by (\ref{2.3}), H\"{o}lder's inequality and Young's inequality, it follows that
 \bes
K_1&=&\int_{\mathbb{R}^3}|u||\nabla u||\nabla_h\nabla u|{\rm d}x\nonumber\\[3mm]
&\leq&{\Big\|}\|u\|_{L_3^p}{\Big\|}_{L_{1, 2}^\alpha} {\Big\|}\|\nabla u\|_{L_3^{\frac{2p}{p-2}}}
{\Big\|}_{L_{1, 2}^{\frac{2\alpha}{\alpha-2}}}
\|\nabla_h\nabla u\|_{L^2}\nonumber\\[3mm]
&\leq& C{\Big\|}\|u\|_{L_3^p}{\Big\|}_{L_{1, 2}^\alpha}
\|\nabla u\|_{L^2}^{1-\frac{1}{p}-\frac{2}{\alpha}}\|\Delta u\|_{L^2}^{\frac{1}{p}}
\|\nabla_h\nabla u\|_{L^2}^{1+\frac{2}{\alpha}}\nonumber\\[3mm]
&\leq& C{\Big\|}\|u\|_{L_3^p}{\Big\|}_{L_{1, 2}^\alpha}^{\frac{2\alpha}{\alpha-2}}
\|\nabla u\|_{L^2}^{\frac{2(p\alpha-\alpha-2p)}{p(\alpha-2)}}
\|\Delta u\|_{L^2}^{\frac{2\alpha}{p(\alpha-2)}}
+\frac{1}{2}\|\nabla_h\nabla u\|_{L^2}^{2}.
 \nonumber
 \ees
In a similar way, one has
 \bes
K_2&=&\int_{\mathbb{R}^3}|u||\nabla b||\nabla_h\nabla b|{\rm d}x\nonumber\\[3mm]
&\leq&{\Big\|}\|u\|_{L_3^p}{\Big\|}_{L_{1, 2}^\alpha} {\Big\|}
\|\nabla b\|_{L_3^{\frac{2p}{p-2}}}{\Big\|}_{L_{1, 2}^{\frac{2\alpha}{\alpha-2}}}
\|\nabla_h\nabla b\|_{L^2}\nonumber\\[3mm]
&\leq& C{\Big\|}\|u\|_{L_3^p}{\Big\|}_{L_{1, 2}^\alpha}
\|\nabla b\|_{L^2}^{1-\frac{1}{p}-\frac{2}{\alpha}}\|\Delta b\|_{L^2}^{\frac{1}{p}}
\|\nabla_h\nabla b\|_{L^2}^{1+\frac{2}{\alpha}}\nonumber\\[3mm]
&\leq& C{\Big\|}\|u\|_{L_3^p}{\Big\|}_{L_{1, 2}^\alpha}^{\frac{2\alpha}{\alpha-2}}
\|\nabla b\|_{L^2}^{\frac{2(p\alpha-\alpha-2p)}{p(\alpha-2)}}\|\Delta b\|_{L^2}^{\frac{2\alpha}{p(\alpha-2)}}
+\frac{1}{2}\|\nabla_h\nabla b\|_{L^2}^{2}.
 \nonumber
 \ees
Collecting the above two  estimates, we have
 \bes
&&\hspace{-1cm}\frac{d}{dt}{\Big(}\|\nabla_h u(t)\|_{L^2}^2+\|\nabla_h b(t)\|_{L^2}^2{\Big)}
+{\Big(}\|\nabla_h \nabla u\|_{L^2}^2
+\|\nabla_h\nabla b\|_{L^2}^2{\Big)}\nonumber\\[3mm]
&\leq&C{\Big\|}\|u\|_{L_3^p}{\Big\|}_{L_{1, 2}^\alpha}^{\frac{2\alpha}{\alpha-2}}
{\Big (}\|\nabla u\|_{L^2}^{\frac{2(p\alpha-\alpha-2p)}{p(\alpha-2)}}
+\|\nabla b\|_{L^2}^{\frac{2(p\alpha-\alpha-2p)}{p(\alpha-2)}}{\Big )}
{\Big (}\|\Delta u\|_{L^2}^{\frac{2\alpha}{p(\alpha-2)}}+\|\Delta b\|_{L^2}^{\frac{2\alpha}{p(\alpha-2)}}{\Big )}.
 \nonumber
 \ees
Integrating in time and using H\"{o}lder's inequality lead to
 \bes
 \begin{array}[b]{ll}
&\hspace{-0.5cm}\D \|\nabla_h u(t)\|_{L^2}^2+\|\nabla_h b(t)\|_{L^2}^2
+\int_0^t{\Big (}\|\nabla_h \nabla u(\tau)\|_{L^2}^2+\|\nabla_h\nabla b(\tau)\|_{L^2}^2{\Big)}{\rm d}\tau\\[3mm]
&\leq\D C{\Big [}\int_0^t{\Big\|}\|u(\tau)\|_{L_3^p}{\Big\|}_{L_{1, 2}^\alpha}^{\frac{2p\alpha}{p\alpha-2p-2}}
{\Big (}\|\nabla u(\tau)\|_{L^2}^{2}
+\|\nabla b(\tau)\|_{L^2}^2{\Big )}{\rm d}\tau{\Big ]}^{\frac{p\alpha-2p-2}{p(\alpha-2)}}\\[3mm]
&\ \ \ \times\D {\Big (}\int_0^t(\|\Delta u(\tau)\|_{L^2}^2
+\|\Delta b(\tau)\|_{L^2}^2){\rm d}\tau{\Big )}^{\frac{\alpha}{p(\alpha-2)}}+\|\nabla_h u_0\|_{L^2}^2+\|\nabla_h b_0\|_{L^2}^2\\[3mm]
&=:CG(t),
 \lbl{2.4}
 \end{array}
 \ees
where
 \bes
 \begin{array}[b]{ll}
G(t)&=\D {\Big [}\int_0^t{\Big\|}\|u(\tau)\|_{L_3^p}{\Big\|}_{L_{1, 2}^\alpha}^{\frac{2p\alpha}{p\alpha-2p-2}}
{\Big (}\|\nabla u(\tau)\|_{L^2}^{2}
+\|\nabla b(\tau)\|_{L^2}^2{\Big )}{\rm d}\tau{\Big ]}^{\frac{p\alpha-2p-2}{p(\alpha-2)}}\\[3mm]
&\ \ \ \times\D {\Big (}\int_0^t(\|\Delta u(\tau)\|_{L^2}^2
+\|\Delta b(\tau)\|_{L^2}^2){\rm d}\tau{\Big )}^{\frac{\alpha}{p(\alpha-2)}}+\|\nabla_h u_0\|_{L^2}^2+\|\nabla_h b_0\|_{L^2}^2.
 \lbl{2.25}
 \end{array}
 \ees

{\em Step 2. Estimates for $\|\nabla u\|_{L^2}$ and $\|\nabla b\|_{L^2}$}.
Taking the inner product of the first and second equations
in (\ref{1.1}) with $-\Delta u$ and $-\Delta b$, respectively, we have
 \bes
 \begin{array}[b]{ll}
&\hspace{-0.5cm} \D \frac{1}{2}\frac{d}{dt}{\Big(}\|\nabla u(t)\|_{L^2}^2
+\|\nabla b(t)\|_{L^2}^2{\Big)}+{\Big(}\|\Delta u\|_{L^2}^2+\|\Delta b\|_{L^2}^2{\Big)}\\[3mm]
&=\D \int_{\mathbb{R}^3}(u\cdot\nabla)u\cdot \Delta u{\rm d}x
-\int_{\mathbb{R}^3}(b\cdot\nabla)b\cdot \Delta u{\rm d}x\\[3mm]
&\ \ \ +\D \int_{\mathbb{R}^3}(u\cdot\nabla)b\cdot \Delta b{\rm d}x
-\int_{\mathbb{R}^3}(b\cdot\nabla)u\cdot \Delta b{\rm d}x\\[3mm]
&=:\D J_1+J_2+J_3+J_4.
 \lbl{2.5}
 \end{array}
 \ees
In the following, we establish the boun{\rm d}s of $J_1$-$J_4$. For the first term $J_1$, we have
 \bess
 \begin{array}{ll}
J_1&=\D \int_{\mathbb{R}^3}(u\cdot \nabla)u\cdot \Delta_h u {\rm d}x
+\int_{\mathbb{R}^3}(u\cdot \nabla)u\cdot \partial_3^2 u {\rm d}x\\[3mm]
&=:J_{11}+J_{12}.
 \end{array}
 \eess
Recalling (\ref{2.20}), one gets
$$J_{11}=\int_{\mathbb{R}^3}(u\cdot \nabla)u\cdot \Delta_h u {\rm d}x\leq
 \int_{\mathbb{R}^3}|u||\nabla u||\nabla_h\nabla u|{\rm d}x.$$
For $J_{12}$, applying integration by parts and $\nabla\cdot u=0$ yield that
 \bes
J_{12}&=&\int_{\mathbb{R}^3}(u\cdot \nabla)u\cdot \partial_3^2 u {\rm d}x=
-\int_{\mathbb{R}^3}(\partial_3 u\cdot \nabla)u\cdot \partial_3 u {\rm d}x\nonumber\\[3mm]
&=&-\sum_{k=1}^2\int_{\mathbb{R}^3}\partial_3u_k\partial_ku\cdot\partial_3u{\rm d}x
-\int_{\mathbb{R}^3}\partial_3u_3\partial_3u\cdot\partial_3u{\rm d}x\nonumber\\[3mm]
&=&-\sum_{k=1}^2\int_{\mathbb{R}^3}\partial_3u_k\partial_ku\cdot\partial_3u{\rm d}x
+\sum_{k=1}^2\int_{\mathbb{R}^3}\partial_ku_k\partial_3u\cdot\partial_3u{\rm d}x\nonumber\\[3mm]
&\leq& \int_{\mathbb{R}^3}|\nabla_h u||\nabla u|^2{\rm d}x.
 \nonumber
 \ees
As a consequence,
$$J_1\leq J_{11}+J_{12}\leq \int_{\mathbb{R}^3}|u||\nabla u||\nabla_h\nabla u|{\rm d}x
+\int_{\mathbb{R}^3}|\nabla_h u||\nabla u|^2{\rm d}x.$$
Similarly,  by (\ref{2.22}), one gets
\bes
J_3&=&\int_{\mathbb{R}^3}(u\cdot \nabla)b\cdot\Delta b{\rm d}x
=\int_{\mathbb{R}^3}(u\cdot \nabla)b\cdot \Delta_h b {\rm d}x
+\int_{\mathbb{R}^3}(u\cdot \nabla)b\cdot \partial_3^2 b {\rm d}x\nonumber\\[3mm]
&\leq&\int_{\mathbb{R}^3}|u||\nabla b||\nabla_h\nabla b|{\rm d}x
-\int_{\mathbb{R}^3}(\partial_3 u\cdot \nabla)b\cdot \partial_3 b {\rm d}x\nonumber\\[3mm]
&=&\int_{\mathbb{R}^3}|u||\nabla b||\nabla_h\nabla b|{\rm d}x
-\sum_{k=1}^2\int_{\mathbb{R}^3}\partial_3u_k\partial_kb\cdot\partial_3b{\rm d}x
+\sum_{k=1}^2\int_{\mathbb{R}^3}\partial_ku_k\partial_3b\cdot\partial_3b{\rm d}x\nonumber\\[3mm]
&\leq&\int_{\mathbb{R}^3}|u||\nabla b||\nabla_h\nabla b|{\rm d}x
+\int_{\mathbb{R}^3}|\nabla u||\nabla b||\nabla_h b|{\rm d}x
+ \int_{\mathbb{R}^3}|\nabla_h u||\nabla b|^2{\rm d}x.
 \nonumber
 \ees
Thanks to (\ref{2.21}), it hol{\rm d}s that
 \bess
 \begin{array}{ll}
J_2+J_4&=\D -\int_{\mathbb{R}^3}(b\cdot \nabla)b\cdot\Delta u{\rm d}x
-\int_{\mathbb{R}^3}(b\cdot \nabla)u\cdot\Delta b{\rm d}x\\[3mm]
&=\D -\int_{\mathbb{R}^3}(b\cdot \nabla)b\cdot\Delta_h u{\rm d}x
-\int_{\mathbb{R}^3}(b\cdot \nabla)u\cdot\Delta_h b{\rm d}x\\[3mm]
&\ \ \ -\D  \int_{\mathbb{R}^3}(b\cdot \nabla)b\cdot\partial_3^2 u{\rm d}x
-\int_{\mathbb{R}^3}(b\cdot \nabla)u\cdot\partial_3^2 b{\rm d}x\\[3mm]
&\leq\D\int_{\mathbb{R}^3}|u||\nabla b||\nabla_h\nabla b|{\rm d}x
+\int_{\mathbb{R}^3}(\partial_3 b\cdot \nabla)b\cdot \partial_3 u {\rm d}x
+\int_{\mathbb{R}^3}(\partial_3 b\cdot \nabla)u\cdot \partial_3 b {\rm d}x\\[3mm]
&\leq\D \int_{\mathbb{R}^3}|u||\nabla b||\nabla_h\nabla b|{\rm d}x
+\sum_{k=1}^2\int_{\mathbb{R}^3}\partial_3 b_k \partial_kb\cdot \partial_3 u {\rm d}x
+\int_{\mathbb{R}^3}\partial_3 b_3\partial_3 b\cdot \partial_3 u{\rm d}x\\[3mm]
&\ \ \ +\D\sum_{k=1}^2\int_{\mathbb{R}^3}\partial_3 b_k \partial_ku\cdot \partial_3 b{\rm d}x
+\int_{\mathbb{R}^3}\partial_3 b_3\partial_3 u\cdot \partial_3 b{\rm d}x\\[3mm]
&=\D\int_{\mathbb{R}^3}|u||\nabla b||\nabla_h\nabla b|{\rm d}x
+\sum_{k=1}^2\int_{\mathbb{R}^3}\partial_3 b_k \partial_kb\cdot \partial_3 u {\rm d}x
-\sum_{k=1}^2\int_{\mathbb{R}^3}\partial_k b_k\partial_3 b\cdot \partial_3 u{\rm d}x\\[3mm]
&\ \ \ +\D\sum_{k=1}^2\int_{\mathbb{R}^3}\partial_3 b_k \partial_ku\cdot \partial_3 b{\rm d}x
-\sum_{k=1}^2\int_{\mathbb{R}^3}\partial_k b_k\partial_3 u\cdot \partial_3 b{\rm d}x\\[3mm]
&\leq\D \int_{\mathbb{R}^3}|u||\nabla b||\nabla_h\nabla b|{\rm d}x
+\int_{\mathbb{R}^3}|\nabla u||\nabla b||\nabla_h b|{\rm d}x
+\int_{\mathbb{R}^3}|\nabla b|^2|\nabla_h u|{\rm d}x.
 \end{array}
 \eess
Plugging the above estimates into (\ref{2.5}) yiel{\rm d}s
 \bes
 \begin{array}[b]{ll}
&\hspace{-0.5cm}\D \frac{1}{2}\frac{d}{dt}{\Big(}\|\nabla u(t)\|_{L^2}^2
+\|\nabla b(t)\|_{L^2}^2{\Big)}+{\Big(}\|\Delta u\|_{L^2}^2
+\|\Delta b\|_{L^2}^2{\Big)}\\[3mm]
&\leq\D \int_{\mathbb{R}^3}|u||\nabla u||\nabla_h\nabla u|{\rm d}x
+\int_{\mathbb{R}^3}|u||\nabla b||\nabla_h\nabla b|{\rm d}x
+\int_{\mathbb{R}^3}|\nabla u|^2|\nabla_h u|{\rm d}x\\[3mm]
&\ \ \ +\D \int_{\mathbb{R}^3}|\nabla u||\nabla b||\nabla_h b|{\rm d}x
+\int_{\mathbb{R}^3}|\nabla b|^2|\nabla_h u|{\rm d}x\\[3mm]
&=:\D H_1+H_2+H_3+H_4+H_5.
 \lbl{2.6}
 \end{array}
 \ees
We now estimate the terms $H_1$-$H_5$ one by one. To estimate the first term $H_1$,
we apply H\"{o}lder's inequality and Young's inequality to derive
 \bes
H_1&=&\int_{\mathbb{R}^3}|u||\nabla u||\nabla_h\nabla u|{\rm d}x\nonumber\\[3mm]
&\leq&{\Big\|}\|u\|_{L_3^p}{\Big\|}_{L_{1, 2}^\alpha} {\Big\|}\|\nabla u\|_{L_3^{\frac{2p}{p-2}}}
{\Big\|}_{L_{1, 2}^{\frac{2\alpha}{\alpha-2}}}
\|\nabla_h\nabla u\|_{L^2}\nonumber\\[3mm]
&\leq& C{\Big\|}\|u\|_{L_3^p}{\Big\|}_{L_{1, 2}^\alpha}
\|\nabla u\|_{L^2}^{1-\frac{1}{p}-\frac{2}{\alpha}}\|\Delta u\|_{L^2}^{1+\frac{2}{\alpha}+\frac{1}{p}}\nonumber\\[3mm]
&\leq& C{\Big\|}\|u\|_{L_3^p}{\Big\|}_{L_{1, 2}^\alpha}^{\frac{2p\alpha}{p\alpha-2p-\alpha}}
\|\nabla u\|_{L^2}^2
+\frac{1}{2}\|\Delta u\|_{L^2}^{2}.
 \nonumber
 \ees
Similarly, one has
 \bes
H_2&=&\int_{\mathbb{R}^3}|u||\nabla b||\nabla_h\nabla b|{\rm d}x\nonumber\\[3mm]
&\leq& C{\Big\|}\|u\|_{L_3^p}{\Big\|}_{L_{1, 2}^\alpha}
\|\nabla b\|_{L^2}^{1-\frac{1}{p}-\frac{2}{\alpha}}\|\Delta b\|_{L^2}^{1+\frac{2}{\alpha}+\frac{1}{p}}\nonumber\\[3mm]
&\leq& C{\Big\|}\|u\|_{L_3^p}{\Big\|}_{L_{1, 2}^\alpha}^{\frac{2p\alpha}{p\alpha-2p-\alpha}}
\|\nabla b\|_{L^2}^2
+\frac{1}{2}\|\Delta b\|_{L^2}^{2}.
 \nonumber
 \ees
On the other hand, the following relation
 \bes
 \|\nabla u\|_{L^4(\mathbb{R}^3)}\leq C\|\nabla u\|_{L^2}^{\frac{1}{4}}
 \|\partial_1\nabla u\|_{L^2}^{\frac{1}{4}}
 \|\partial_2\nabla u\|_{L^2}^{\frac{1}{4}}
 \|\partial_3\nabla u\|_{L^2}^{\frac{1}{4}}
 \nonumber
 \ees
derived from Lemma \ref{lm1.1} together with H\"{o}lder's inequality ensures that
 \bes
H_3&=&\int_{\mathbb{R}^3}|\nabla u|^2|\nabla_h u|{\rm d}x\leq \|\nabla_h u\|_{L^2}\|\nabla u\|_{L^4}^2\nonumber\\[3mm]
&\leq& C\|\nabla_h u\|_{L^2}\|\nabla u\|_{L^2}^{\frac{1}{2}}\|\nabla_h\nabla u\|_{L^2}\|\Delta u\|_{L^2}^{\frac{1}{2}},
 \nonumber
 \ees
 \bes
H_4&=&\int_{\mathbb{R}^3}|\nabla u||\nabla b||\nabla_h b|{\rm d}x
\leq\|\nabla_h b\|_{L^2}\|\nabla u\|_{L^4}\|\nabla b\|_{L^4}\nonumber\\[3mm]
&\leq& C\|\nabla_h b\|_{L^2}\|\nabla u\|_{L^2}^{\frac{1}{4}}
\|\nabla_h\nabla u\|_{L^2}^{\frac{1}{2}}\|\partial_3\nabla u\|_{L^2}^{\frac{1}{4}}
\|\nabla b\|_{L^2}^{\frac{1}{4}}
\|\nabla_h\nabla b\|_{L^2}^{\frac{1}{2}}\|\partial_3\nabla b\|_{L^2}^{\frac{1}{4}}\nonumber\\[3mm]
&\leq& C\|\nabla_h b\|_{L^2}(\|\nabla u\|_{L^2}^{\frac{1}{2}}
+\|\nabla b\|_{L^2}^{\frac{1}{2}})(\|\nabla_h\nabla u\|_{L^2}+\|\nabla_h\nabla b\|_{L^2})
(\|\Delta u\|_{L^2}^{\frac{1}{2}}+\|\Delta b\|_{L^2}^{\frac{1}{2}}),
 \nonumber
 \ees
and
 \bes
H_5&=&\int_{\mathbb{R}^3}|\nabla b|^2|\nabla_h u|{\rm d}x\leq \|\nabla_h u\|_{L^2}\|\nabla b\|_{L^4}^2\nonumber\\[3mm]
&\leq& C\|\nabla_h u\|_{L^2}\|\nabla b\|_{L^2}^{\frac{1}{2}}\|\nabla_h\nabla b\|_{L^2}\|\Delta b\|_{L^2}^{\frac{1}{2}}.
 \nonumber
 \ees
Thus, substituting the above estimates into (\ref{2.6}), we get
 \bess
 \begin{array}{ll}
&\hspace{-0.5cm}\D \frac{d}{dt}{\Big(}\|\nabla u(t)\|_{L^2}^2
+\|\nabla b(t)\|_{L^2}^2{\Big)}+{\Big(}\|\Delta u\|_{L^2}^2
+\|\Delta b\|_{L^2}^2{\Big)}\\[3mm]
&\leq\D C(\|\nabla_h u\|_{L^2}+\|\nabla_h b\|_{L^2})(\|\nabla u\|_{L^2}^{\frac{1}{2}}
+\|\nabla b\|_{L^2}^{\frac{1}{2}})(\|\nabla_h\nabla u\|_{L^2}+\|\nabla_h\nabla b\|_{L^2})\\[3mm]
&\ \ \ \times\D (\|\Delta u\|_{L^2}^{\frac{1}{2}}+\|\Delta b\|_{L^2}^{\frac{1}{2}})
+ C{\Big\|}\|u\|_{L_3^p}{\Big\|}_{L_{1, 2}^\alpha}^{\frac{2p\alpha}{p\alpha-2p-\alpha}}
(\|\nabla u\|_{L^2}^2+\|\nabla b\|_{L^2}^2).
 \end{array}
 \eess
Integrating in time and using H\"{o}lder's inequality, we have
 \bes
 \begin{array}[b]{ll}
&\hspace{-0.5cm} \D{\Big(}\|\nabla u(t)\|_{L^2}^2
+\|\nabla b(t)\|_{L^2}^2{\Big)}+\int_0^t{\Big(}\|\Delta u(\tau)\|_{L^2}^2
+\|\Delta b(\tau)\|_{L^2}^2{\Big)}{\rm d}\tau\\[3mm]
&\leq \D  C\int_0^t{\Big\|}\|u(\tau)\|_{L_3^p}{\Big\|}_{L_{1, 2}^\alpha}^{\frac{2p\alpha}{p\alpha-2p-\alpha}}
(\|\nabla u(\tau)\|_{L^2}^2+\|\nabla b(\tau)\|_{L^2}^2){\rm d}\tau\\[3mm]
&\ \ \ + CY(t)+ (\|\nabla u_0\|_{L^2}^2+\|\nabla b_0\|_{L^2}^2),
 \lbl{2.7}
 \end{array}
 \ees
where
 \bess
 \begin{array}{ll}
Y(t)&=\D \int_0^t(\|\nabla_h u(\tau)\|_{L^2}+\|\nabla_h b(\tau)\|_{L^2})
(\|\nabla u(\tau)\|_{L^2}^{\frac{1}{2}}+\|\nabla b(\tau)\|_{L^2}^{\frac{1}{2}})\\[3mm]
&\ \ \ \times\D(\|\nabla_h\nabla u(\tau)\|_{L^2}+\|\nabla_h\nabla b(\tau)\|_{L^2})
(\|\Delta u(\tau)\|_{L^2}^{\frac{1}{2}}+\|\Delta b(\tau)\|_{L^2}^{\frac{1}{2}}){\rm d}\tau.
 \end{array}
 \eess
Next, we establish the bound of $Y(t)$.
 \bes
 \begin{array}[b]{ll}
Y(t)&\leq\D (\|\nabla_h u\|_{L_t^\infty L^2}+\|\nabla_h b\|_{L_t^\infty L^2})
(\|\nabla u\|_{L_t^2L^2}^{\frac{1}{2}}+\|\nabla b\|_{L_t^2L^2}^{\frac{1}{2}})\\[3mm]
&\ \ \ \times\D(\|\nabla_h\nabla u\|_{L_t^2L^2}+\|\nabla_h\nabla b\|_{L_t^2L^2})
(\|\Delta u\|_{L_t^2L^2}^{\frac{1}{2}}+\|\Delta b\|_{L_t^2L^2}^{\frac{1}{2}})\\[3mm]
&\leq\D CG(t)(\|\Delta u\|_{L_t^2L^2}^{\frac{1}{2}}+\|\Delta b\|_{L_t^2L^2}^{\frac{1}{2}})\\[3mm]
&\leq\D CG(t){\Big (}\int_0^t(\|\Delta u(\tau)\|_{L^2}^2+\|\Delta b(\tau)\|_{L^2}^2){\rm d}\tau{\Big)}^{\frac{1}{4}},
 \lbl{2.8}
 \end{array}
 \ees
where we have used (\ref{2.0}) and (\ref{2.4}) and $G(t)$ is defined in (\ref{2.25}). Substituting (\ref{2.25}) into (\ref{2.8}) and
applying H\"{o}lder's inequality and Young's inequality again, we deduce
 \bes
 \begin{array}[b]{ll}
Y(t)&\leq\D C(\|\nabla_hu_0\|_{L^2}^2+\|\nabla_hb_0\|_{L^2}^2){\Big (}\int_0^t(\|\Delta u(\tau)\|_{L^2}^2+\|\Delta b(\tau)\|_{L^2}^2){\rm d}\tau{\Big)}^{\frac{1}{4}}\\[3mm]
&\ \ \ +\D C{\Big (}\int_0^t{\Big\|}\|u(\tau)\|_{L_3^p}{\Big\|}_{L_{1, 2}^\alpha}^{\frac{2p\alpha}{p\alpha-2p-\alpha}}
(\|\nabla u(\tau)\|_{L^2}^2+\|\nabla b(\tau)\|_{L^2}^2){\rm d}\tau{\Big )}^{\frac{p\alpha-2p-\alpha}{p(\alpha-2)}}\\[3mm]
&\ \ \ \times\D {\Big (}\int_0^t(\|\Delta u(\tau)\|_{L^2}^2
+\|\Delta b(\tau)\|_{L^2}^2){\rm d}\tau{\Big)}^{\frac{\alpha}{p(\alpha-2)}+\frac{1}{4}}\\[3mm]
&\leq\D  C(\|\nabla u_0\|_{L^2}^{\frac{8}{3}}+\|\nabla b_0\|_{L^2}^{\frac{8}{3}})+\frac{1}{4}\int_0^t(\|\Delta u(\tau)\|_{L^2}^2+\|\Delta b(\tau)\|_{L^2}^2){\rm d}\tau\nonumber\\
&\ \ \ +\D C{\Big (}\int_0^t{\Big\|}\|u(\tau)\|_{L_3^p}{\Big\|}_{L_{1, 2}^\alpha}^{\frac{8p\alpha}{3p\alpha-6p-4\alpha}}
(\|\nabla u(\tau)\|_{L^2}^2+\|\nabla b(\tau)\|_{L^2}^2){\rm d}\tau{\Big )}^{\frac{3p\alpha-6p-4\alpha}{4p(\alpha-2)}}\\[3mm]
&\ \ \ \times\D {\Big (}\int_0^t(\|\nabla u(\tau)\|_{L^2}^2
+\|\nabla b(\tau)\|_{L^2}^2){\rm d}\tau{\Big )}^{\frac{1}{4}}\\[3mm]
&\ \ \ \times\D {\Big (}\int_0^t(\|\Delta u(\tau)\|_{L^2}^2+\|\Delta b(\tau)\|_{L^2}^2){\rm d}\tau{\Big)}^{\frac{\alpha}{p(\alpha-2)}+\frac{1}{4}}\nonumber\\[3mm]
&\leq\D  C(\|\nabla u_0\|_{L^2}^{\frac{8}{3}}+\|\nabla b_0\|_{L^2}^{\frac{8}{3}})+\frac{1}{2}\int_0^t(\|\Delta u(\tau)\|_{L^2}^2+\|\Delta b(\tau)\|_{L^2}^2){\rm d}\tau\\[3mm]
&\ \ \ +\D C\int_0^t{\Big\|}\|u(\tau)\|_{L_3^p}{\Big\|}_{L_{1, 2}^\alpha}^{\frac{8p\alpha}{3p\alpha-6p-4\alpha}}
(\|\nabla u(\tau)\|_{L^2}^2+\|\nabla b(\tau)\|_{L^2}^2){\rm d}\tau,
 \lbl{2.9}
 \end{array}
 \ees
where we have used (\ref{2.0}) again. Finally, by (\ref{2.7}) and (\ref{2.9}), we have
 \bes
&&\hspace{-1cm}{\Big(}\|\nabla u(t)\|_{L^2}^2
+\|\nabla b(t)\|_{L^2}^2{\Big)}+\frac{1}{2}\int_0^t{\Big(}\|\Delta u(\tau)\|_{L^2}^2
+\|\Delta b(\tau)\|_{L^2}^2{\Big)}{\rm d}\tau\nonumber\\[3mm]
&\leq&C(1+\|\nabla u_0\|_{L^2}^{\frac{8}{3}}+\|\nabla b_0\|_{L^2}^{\frac{8}{3}})
+C\int_0^tF(\tau)(\|\nabla u(\tau)\|_{L^2}^2+\|\nabla b(\tau)\|_{L^2}^2){\rm d}\tau,
 \ees
here $F(\tau)={\Big \|}\|u(\tau)\|_{L_3^p}{\Big\|}_{L_{1, 2}^\alpha}^{\frac{2p\alpha}{p\alpha-2p-\alpha}}
+{\Big \|}\|u(\tau)\|_{L_3^p}{\Big\|}_{L_{1, 2}^\alpha}^{\frac{8p\alpha}{3p\alpha-6p-4\alpha}}$.
Thus, Gronwall's inequality guarantees that
 \bes
&&\hspace{-1cm}{\Big(}\|\nabla u(t)\|_{L^2}^2
+\|\nabla b(t)\|_{L^2}^2{\Big)}+\int_0^t{\Big(}\|\Delta u(\tau)\|_{L^2}^2+\|\Delta b(\tau)\|_{L^2}^2{\Big)}{\rm d}\tau\nonumber\\[3mm]
&\leq& C{\Big(}e^{C\int_0^tF(\tau){\rm d}\tau}+t{\Big )}
 \nonumber
 \ees
for all $t\in [0, T^\ast)$.  By means of  condition (\ref{2.30}),
it follows that the $H^1$ norm of the strong solution $(u, b)$ is bounded on the maximal interval of existence $[0, T^\ast)$.
This en{\rm d}s the proof of Theorem \ref{th1.1}.

\subsection{Proof of Corollary \ref{c1}}
In this subsection, we will prove Corollary \ref{c1}. Firstly, when
$\gamma\in [0,\frac{1}{\alpha})$, the desired result follows
directly from the embedding theorem
$${\Big\|}\|u\|_{L_i^{\frac{\alpha}{1-\alpha\gamma}}}{\Big\|}_{L_{j, k}^{\alpha}}\leq C\|\Lambda_i^\gamma u\|_{L^\alpha}.$$
When $\gamma\in (\frac{1}{\alpha}, 1]$, our objective  is to prove that
$$\displaystyle\int_0^t\|\Lambda_i^\gamma u(\tau)\|_{L^\alpha}^\beta {\rm d}\tau<\infty,
~~~\frac{2}{\beta}+\frac{3}{\alpha}\leq \frac{3}{4}
+\frac{3}{2\alpha},~~\alpha\in (2,\infty)$$
is a sufficient condition.
It is easy to check that the integral term
$\int_0^t\|\|u(\tau)\|_{L_i^\infty}\|_{L_{j, k}^{\delta}}^{\frac{8\delta}{3(\delta-2)}}{\rm d}\tau$
satisfies the condition of Theorem \ref{th1.1} with $\delta\in(2,\infty)$.
By Lemma \ref{lm1.1} and
the interpolation theorem,
we obtain that for $\delta\in [\frac{(2\gamma+1)\alpha-2}{\gamma\alpha}, \alpha]$,
 \bes
{\Big\|}\|u\|_{L_i^\infty}{\Big\|}_{L_{j,k}^\delta}
&\leq&C{\Big\|}\|u\|_{L_i^2}^\theta\|\Lambda_i^\gamma u\|_{L_i^\alpha}^{1-\theta}{\Big\|}_{L_{j,k}^\delta}\nonumber\\[3mm]
&\leq& C{\Big\|}\|u\|_{L_i^2}^\theta{\Big\|}_{L_{j,k}^{\frac{p}{\theta}}}{\Big\|}\|\Lambda_i^\gamma u\|_{L_i^\alpha}^{1-\theta}{\Big\|}_{L_{j,k}^{\frac{\alpha}{1-\theta}}}\nonumber\\[3mm]
&=& C{\Big\|}\|u\|_{L_i^2}{\Big\|}_{L_{j,k}^{p}}^\theta{\Big\|}\|\Lambda_i^\gamma u\|_{L_i^\alpha}{\Big\|}_{L_{j,k}^{{\alpha}}}^{1-\theta}\nonumber\\[3mm]
&\leq& C{\Big\|}\|u\|_{L_{j,k}^{p}}{\Big\|}_{L_i^2}^\theta\|\Lambda_i^\gamma u\|_{L^\alpha}^{1-\theta}\nonumber\\[3mm]
&\leq& C{\Big\|}\|u\|_{L_{j,k}^{2}}^{\frac{2}{p}}\|\partial_ju\|_{L_{j,k}^{2}}^{\frac{p-2}{2p}}
\|\partial_ku\|_{L_{j,k}^{2}}^{\frac{p-2}{2p}}{\Big\|}_{L_i^2}^\theta\|\Lambda_i^\gamma u\|_{L^\alpha}^{1-\theta}\nonumber\\[3mm]
&\leq& C\|u\|_{L^2}^{\frac{2\theta}{p}}\|\nabla u\|_{L^2}^{\frac{(p-2)\theta}{p}}\|\Lambda_i^\gamma u\|_{L^\alpha}^{1-\theta},
 \nonumber
 \ees
where $\frac{1}{\delta}=\frac{\theta}{p}+\frac{1-\theta}{\alpha}$ with $\theta=\frac{2(\gamma\alpha-1)}{2(\gamma\alpha-1)+\alpha}$.
In the forth line we have used the fact $$\delta\geq\frac{(2\gamma+1)\alpha-2}{\gamma\alpha}\Leftrightarrow p\geq 2.$$
Therefore, H\"{o}lder's inequality enables us to deduce that
 \bes
 \int_0^t{\Big\|}\|u(\tau)\|_{L_i^\infty}{\Big\|}_{L_{j, k}^{\delta}}^{\frac{8\delta}{3(\delta-2)}}{\rm d}\tau
&\leq&C\int_0^t\|u(\tau)\|_{L^2}^{\frac{16\delta\theta}{3p(\delta-2)}}\|\nabla u(\tau)\|_{L^2}^{\frac{8\delta(p-2)\theta}{3p(\delta-2)}}\|\Lambda_i^\gamma u(\tau)\|_{L^\alpha}^{\frac{8\delta(1-\theta)}{3(\delta-2)}}{\rm d}\tau\nonumber\\[3mm]
&\leq&C\|u\|_{L_t^\infty L^2}^{\frac{16\delta\theta}{3p(\delta-2)}}\|\nabla u\|_{L_t^2L^2}^{\frac{8\delta(p-2)\theta}{3p(\delta-2)}}{\Big (}\int_0^t\|\Lambda_i^\gamma u(\tau)\|_{L^\alpha}^{\frac{8\delta(1-\theta)\eta}{3(\delta-2)}}{\rm d}\tau{\Big )}^{\frac{1}{\eta}},
 \nonumber
 \ees
where $\eta=\frac{3p(\delta-2)}{3p(\delta-2)-4\delta(p-2)\theta}$.

According to the relation $\frac{1}{\delta}=\frac{\theta}{p}+\frac{1-\theta}{\alpha}$ and $\theta=\frac{2(\gamma\alpha-1)}{2(\gamma\alpha-1)+\alpha}$, we get
 \bes
 \frac{8\delta(1-\theta)\eta}{3(\delta-2)}&=&\frac{8p\delta(1-\theta)}{3p(\delta-2)-4\delta(p-2)\theta}\nonumber\\[3mm]
&=&\frac{8\delta(1-\theta)}{3(\delta-2)-4\delta(1-\frac{2}{p})\theta}\nonumber\\[3mm]
&=&\frac{8\alpha(1-\theta)}{3\alpha+\frac{2\alpha}{\delta}-4\theta\alpha-8(1-\theta)}\nonumber\\[3mm]
&=:&H(\delta).
 \nonumber
 \ees
Since the value of $\theta$ is independent of $\delta$,   it is clear that $H(\delta)$ is a strictly increasing function in terms of $\delta$ on interval $[2,\alpha]$. This together with the fact $H(\frac{(2\gamma+1)\alpha-2}{\gamma\alpha})=\frac{8\alpha}{3\alpha-6}$ and Theorem \ref{th1.1} gives the desired result of Corollary \ref{c1}.


\begin{thebibliography}{}
 \bibliographystyle{siam}
 \setlength{\baselineskip}{12pt}{\small



 \bibitem{L1}
A. Lifschitz, Magnetohydrodynamics and Spectral Theory, Developments in Electromagnetic Theory and Applications, 4, Kluwer Academic Publishers Group, Dordrecht, 1989.

 \bibitem{ST}
M. Sermange, R. Temam,  Some mathematical questions related to the MHD equations, Commun. Pure  Appl. Math. 36 (1983)  635-664.

 \bibitem{CMZ}
Q. Chen, C. Miao, Z. Zhang, On the regularity criterion of weak solution for the 3D viscous magnetohydrodynamics
equations, Commun. Math. Phys.  284 (2008) 919-930.

 \bibitem{CW}
C. Cao, J. Wu, Two regularity criteria for the 3D MHD equations, J. Differential Equations  248 (2010) 2263-2274.

 \bibitem{G}
S. Gala, Extension criterion on regularity for weak solutions to the 3D MHD equations, Math. Meth. Appl.
Sci.   33  (2010) 1496-1503.

 \bibitem{HW1}
C. He, Y. Wang,  On the regularity for weak solutions to the magnetohydrodynamic equations, J. Differential Equations  238  (2007)  1-17.

 \bibitem{HW2}
C. He, Y. Wang, Remark on the regularity for weak solutions to the magnetohydrodynamic equations,
Math. Meth. Appl. Sci.  31   (2008)  1667-1684.

 \bibitem{HX1}
C. He, Z. Xin, On the regularity of solutions to the magnetohydrodynamic equations, J. Differential Equations 213  (2005) 235-254.

 \bibitem{HX2}
C. He, Z. Xin, Partial regularity of suitable weak solutions to the incompressible magnetohydrodynamic equations, J. Funct. Anal. 227 (2005) 113-152.

  \bibitem{JL}
E. Ji, J. Lee, Some regularity criteria for the 3D incompressible magnetohydrodynamics, J. Math. Anal. Appl.  369 (2010) 317-322.

\bibitem{W1}
J. Wu, Viscous and inviscid magnetohydrodynamics equations, J. Anal. Math. 73 (1997) 251-265.

 \bibitem{W2}
J. Wu, Boun{\rm d}s and new approaches for the 3D MHD equations, J. Nonlinear Sci. 12 (2002) 395-413.

 \bibitem{W3}
J. Wu, Regularity results for weak solutions of the 3D MHD equations, Discrete Contin. Dyn. Syst. A 10 (2004) 543-556.

 \bibitem{Z1}
Y. Zhou,
Remarks on regularities for the 3D MHD equations, Discrete Contin. Dyn. Syst.  12  (2005) 881-886.

 \bibitem{Z6}
Y. Zhou, Regularities criteria for the 3D MHD equations in term of pressure, Int. J. Non-Linear Mech. 41 (2006) 1174-1180.

\bibitem{JZ}
X. Jia, Y. Zhou, Regularity criteria for the 3D MHD equations involving partial components, Nonlinear Analysis: RWA 13 (2012) 410-418.

 \bibitem{LD}
H. Lin, L. Du, Regularity criteria for incompressible magnetohydrodynamics equations in three dimensions, Nonlinearity 26 (2013) 219-239.


\bibitem{CT1}
C. Cao, E. S. Titi, Regularity criteria for the three-dimensional Navier-Stokes equations, Indiana Univ. Math. J.  57  (2008) 2643-2661.


 \bibitem{CT2}
C. Cao, E. S. Titi, Global regularity criterion for the 3-D Navier-Stokes equations involving one entry of the velocity gradent tensor, Arch. Rational Mech. Anal. 202 (2011) 919-932.



\bibitem{Z}
X. Zheng, A regularity criterion for the tridimensional Navier-Stokes equations in term of one velocity component, J. Differential Equations 256 (2014) 283-309.

\bibitem{A}
R. A. Adams, Sobolev Spaces, Academic Press, New York, 1975.

\bibitem{DL}
G. Duvaut, J.-L. Lions, In\'{e}quations en thermo\'{e}lasticit\'{e} et magn\'{e}tohydrodynamique, Arch. Rational Mech. Anal. 46 (1972) 241-279.


 \bibitem{L}
O. A. Ladyzhenskaya, Mathematical Theory of Viscous Incompressible Flow, 2nd ed., Gordon and Breach, New York, 1969,
English translation.






}

\end{thebibliography}
\end {document}